\documentclass[12pt]{amsart}
\usepackage{amssymb,  xypic, euscript}
\def\cites{\cite} 
\usepackage{times}

\newif\ifdraft \draftfalse
\ifdraft
      \usepackage{showkeys}
\fi
\newif\ifdraft\draftfalse
\newcommand\labl[1]{
  \ifdraft
     \ifmmode\renewcommand{\theequation}{\arabic{equation}-#1}
     \else{\bf(#1)}
     \fi
  \fi
  \label{#1}
}
\newcommand\Aff{\mathbb{A}}
 \newcommand\CC{\mathbb{C}}

\newcommand\NN{\mathbb{N}}  
 \newcommand\QQ{\mathbb{Q}} 
\newcommand\RR{\mathbb{R}}

\newcommand\ZZ{\mathbb{Z}} 
\newtheorem{theorem}{Theorem}[section]
\newtheorem{lemma}[theorem]{Lemma}
\newtheorem{proposition}[theorem]{Proposition}

\newtheorem{definition-proposition}[theorem]{Definition-Proposition}

\theoremstyle{definition}

\newtheorem{definition}[theorem]{Definition}

\theoremstyle{remark}
\newtheorem{remark}[theorem]{Remark}

\newcommand\coords{\mathcal{O}}

\newcommand\periods{\mathbf{P}}
\newcommand\Ralg{\RR_{\rm{alg}}}

\begin{document}
\title{Periods and Igusa Zeta functions}
\author{Prakash Belkale}
\address{Mathematics Department\\
CB \#3250, Phillips Hall\\
UNC-Chapel Hill\\
Chapel Hill, NC 27599}
\email{belkale@email.unc.edu}
\author{Patrick Brosnan}
\address{UCLA Mathematics Department\\
Box 951555\\
Los Angeles, CA 90095-1555}
\email{pbrosnan@math.ucla.edu}
\begin{abstract} 
We show that the coefficients in the Laurent series of Igusa
Zeta functions $I(s)=\int_C f^s \omega $ are periods. This will be used to show in a subsequent paper (by P. Brosnan) that certain numbers  occurring in Feynman amplitudes (upto gamma factors) are periods. 
\end{abstract}
\maketitle
\section*[Introduction]{Introduction}
In their paper \cite{KontZagier}, Kontsevich and Zagier give an
elementary definition of a period integral as an absolutely convergent
integral of a rational function over a subset of $\RR^n$ defined by
polynomial inequalities and equalities.  They then show that some of
the most important quantities in mathematics are periods, and sketch a
proof that their notion of a period agrees with the more elaborate
notion that algebraic geometers have studied since Riemann and
Weirstra\ss.  The last chapter links periods to the ``framed motives''
studied by A.~Goncharov and proposes a structure of a
torsor on a certain set of framed motives.  The paper is full of interesting
examples, however, its main purpose seems to be to justify the following:  
\theoremstyle{theorem}
\newtheorem{principle}[theorem]{Philosophical Principle}
\begin{principle}
Whenever you meet a new number, and have
decided (or convinced yourself) that it is transcendental, try to
figure out whether it is a period.
\end{principle}

A very interesting class of numbers arises naturally in quantum field theory which we want to prove to be periods.  Namely, if $I(D)$ is a Feynman amplitude coming from a
scalar field theory corresponding to a Feynman integral with all
parameters in $\QQ$, then $I(D)=G(D)J(D)$ where $G(D)$ is a relatively
simple gamma factor and $J(D)$ is a meromorphic function such that the
coefficients in the Laurent series expansion of $J(D)$ at $D=D_0$ are
periods for $D_0$ any integer.  We remark that this confirms (albeit
in a very weak sense) the fact noticed by Kreimer and Broadhurst that
the principal parts of the Laurent series for primitive diagrams often
have coefficients which are multiple zeta values.  

Once dimensional regularization is understood precisely, the proof of
the period-icity of these numbers follows  from an analogous result for Igusa zeta functions which may be interesting in its own right  and which can be explained directly in terms of pure mathematics. We will approach the purely mathematical side in this paper. A subsequent paper by one of us (Brosnan) will explain
the physics of dimensional regularization and on how the theorem on Igusa Zeta
functions ties up with regularization.

To explain this result on Igusa Zeta function,
let $\Delta_n\subset\Bbb{R}^{n+1}$ be the standard $n$-simplex equipped
with the $n$-form 
$$\omega=dx_1\wedge\cdots\wedge dx_n.
$$  
Let $f\in\Bbb{R}[x_0,\ldots, x_n]$ be a polynomial function which is non-negative
on $\Delta$.   Then, according to results of Atiyah \cite{Atiyah70} and
Bernstein and Gelfand \cite{BernsteinGelfand69}, the function 
$$
I(s)= \int_{\Delta_n} f^s\omega
$$
is meromorphic on the complex $s$-plane with isolated singularities.
These functions are called {\em Igusa zeta functions}.
Our main theorem concerning them is the following:
\begin{theorem}
\label{thm:0.1}
Suppose that $f\in\QQ[x_0,\ldots, x_n]$ is a polynomial with rational coeffiecients and 
let $s_0$ be an integer.   Let 
$$
I(s)=\sum_{i\geq N} a_i (s-s_0)^i
$$
be the Laurent series expansion of $I(s)$ at $s_0$.  Then the $a_i$ are 
periods.
\end{theorem}

The above result suffices to show that most of the numbers 
investigated by Kreimer and Broadhurst are, in fact, periods.
However,  it will be convenient to prove a version of this result which 
is more general in the following two senses:  (a) the simplex $\Delta_n$ 
can be taken to be a general semi-algebraic set defined over 
$\overline{\QQ}$, and (b) the function $f$ can be taken to lie in 
the function field $\overline{\QQ}(x_0,\ldots, x_n)$.   For (b) we will need
to use a more general definition of $I(s)$ than the one in  
\cites{Atiyah70, BernsteinGelfand69}.  However, the generalization is 
necessary to handle many of the Feynman amplitudes with infrared
divergences considered by physicists.

We will use the symbol $\periods$ to denote the $\QQ$-algebra of periods 
and use the definition of a period that appears in  \cite{KontZagier}.  
For the convenience of the reader, we also paraphrase this definition in 
(\ref{def:period}).

We thank A.Goncharav, D.Kreimer, M.V.Nori and  H.Rossi  for useful communication.
The idea of using Picard-Fuchs equations in Theorem 1.8 comes from discussions with Madhav Nori. This idea is `standard' when studying periods of powers of functions, but it came somewhat of a surprise that there were no `Gamma factors' at the end. Also historical precedents to the `functional equation' in theorem 1.8 should be noted. They appear in Bernstein's paper
\section[igusa]{Igusa Zeta Functions}

\subsection{Atiyah's Theorem}  Let $X$ be a smooth complex algebraic
variety defined over $\RR$.  Let $X(\RR)$ denote
the real points of $X$
and let $G$ be a semi-algebraic subset of $X(\RR)$
defined by inequalities 
\begin{equation}
  \label{eq:G} 
G=\{x\in X(\RR)\, |\, g_i(x)\geq 0\, \text{for all $i$}\}
\end{equation}
where here the $g_i$ are real-analytic funtions on $X(\RR)$.  
Let $f$ be a real-analytic function on $X(\RR)$ which is non-negative
and not identically zero on $G$.  Let $\Gamma$ denote the characteristic
function of $G$.
In this notation Atiyah's theorem~\cite{Atiyah70} 
can be stated as follows.
\begin{theorem}[Atiyah]
  \label{thm:Atiyah}
  The function $f^s\Gamma$, which is locally integrable
  for $\Re (s)>0$, extends analytically to a distribution on $X$ which
  is a meromorphic function of $s$ in the whole complex plane.  Over
  any relatively compact open set $U$ in $X$ the poles of $f^s\Gamma$
  occur at the points of the form $-r/N, r=1,2,\cdots$, where $N$ is a
  fixed integer (depending on $f$ and $U$) and the order of any pole
  does not exceed the dimension of $X$.  Moreover, $f^0\Gamma =
  \Gamma$.
\end{theorem}

\subsection{Sem-algebraic Sets}
The following definition is given in~\cite{bochnak}.
\begin{definition}
  \label{def:semalg}
A region $C\subset \RR^n$ is {\em semi-algebraic} if it is a union of 
intersections of sets of the form $\{x\in X(\RR) | f(x)>0\}$ or 
$\{x\in X(\RR) | f(x)=0\}$ with $f\in\RR[x_1,\ldots, x_n]$.  
\end{definition}

We will say that $C\subset \RR^n$ is {\em semi-arithmetic} if the 
functions $f$ appearing in the definition are 
in $\Ralg [X_1,\ldots, X_n]$ with $\Ralg =\RR\cap\overline{\QQ}$.

\begin{definition}\label{def:period}
A period is a number whose real and imaginary parts
are given by absolutely convergent integrals of the form 
$\int_C f d\mu$ where $C\subset\RR^n$ is a semi-arithmetic set, 
$f\in\Ralg(x_1,\ldots, x_n)$ and 
$\mu$ is Lebegue measure on $\RR^n$.
\end{definition}

Assume that $X$ is a variety defined over $\RR$. 
For $f\in\RR[X]$,
let $X_{f\geq 0}$ denote the set 
$$
\{x\in X(\RR) | f(x)\geq 0 \}.
$$  
Every point $x\in X$ has an affine neighborhood $V$ which is
isomorphic to a closed subset of $\Aff^n$.  Following \cite{bochnak},
we say that a set $C\subset X(\RR)$ is {\em semi-algebraic} if $C\cap
V$ (considered as a subset of $\RR^n$) is semi-algebraic for every
such affine neighborhood $V$.  If $X$ and $x$ are defined over $\Ralg$, we 
can find a $V$ also defined over $\Ralg$.  We say that
$C$ is {\em semi-arithmetic} if $C\cap V$ is semi-arithmetic for all such $V$.
Clearly, $X_{f\geq 0}$ is semi-algebraic for $f\in\RR[X]$ and
semi-arithmetic for $f\in\Ralg [X]$.

Let $C\subset X(\Bbb{R})$ be a semi-algebraic (resp. semi-arithmetic) set
contained in a dimension $n$ variety $X$, which contains an open (in the usual topology) subset of $X(\Bbb{R})$. It is known that the interior of $C$ contains a semi-algebraic
(resp. semi-arithmetic) dense open subset $U\subset C$ which is smooth
and orientable.  (This follows from Proposition 2.9.10 of~\cite{bochnak}.) 
By a {\em pre-orientation} of $C$, we mean a
choice of such a subset $U$ along with an orientation of $U$.  If
$\omega\in\Omega^n(X)$ is a differential form and $C$ is pre-oriented,
then we make the definition 

\begin{equation}
  \label{eq:preor}
  \int_C \omega \stackrel{\rm{def}}{=} \int_U \omega.
\end{equation}

If $C\subset\RR^n$ then the interior of $C$ is smooth and comes with a 
canonical pre-orientation inherited from the standard orientation on $\RR^n$.
For $C$ compact, the orientation gives a class in 
$\sigma\in H_n(C, \partial C)$ 
where $\partial C$ is the topological boundary of $C$.  
To use Atiyah's theorem in the context of semi-algebraic sets, 
we need to be able to convert an integral $\int_C \omega$ over
an arbitrary semi-algebraic set into a sum of integrals over sets of the 
form of the set $G$ in (\ref{eq:G}). The following lemma is needed to this end.
\begin{lemma} Let $f_i$ ($1\leq i \leq n$) and $g_j$ 
($1\leq j\leq m$) 
be two sets of functions in $\Bbb{R} [X]$.  Let $U=U_1\cup U_2$ be an 
oriented open set with 
\begin{eqnarray}
U_1 &=& \{x\in X | f_i>0\, 1\leq i\leq n \}, \\
U_2 &=& \{x\in X | g_i>0\, 1\leq j\leq m \}.
\end{eqnarray}
Consider strings of the form
$${\bf{e}}=(a_1,\dots,a_n,b_1,\dots,b_m)$$

where the $a_i,b_j$ are in $\{+1,-1\}$ and either all the $a$'s are $+1$
or all the $b$'s are $+1$.

Consider

$$U_{\bf{e}}=\{x\in X | a_i f_i>0, b_jg_j>0, 1\leq i\leq n;1\leq j\leq m \}$$

Then, for a form $\omega\in\Omega^n (X)$ (with $n=\dim X$),
\begin{equation}
\int_U \omega = \sum_{{\bf{e}}}\int_{U_{\bf{e}}} \omega 
\end{equation}
where the ${\bf{e}}$ are subject to the above constraints.
\end{lemma}

Since our domains of integration are going to be semi-algebraic sets, we need a more flexible definition of Periods. This is equivalent to the definition of periods above (~\ref{def:period}), a proof is sketched in \cite{KontZagier}, page 3,31. The definition below is the definition of periods for the purposes of this paper.
\begin{definition}\label{def:period2}
Let $X$ be a smooth algebraic variety of dimension $d$ defined over $\Bbb{Q}$
, $D\subset X$ a divisor with normal crossings, $\omega\in \Omega^d(X)$ an
algebraic differential form on $X$ of the top degree, and $\gamma\in H_d(X(\Bbb{C}),D(\Bbb{C});\Bbb{Q})$ a (homology class of) a singular $d-$ chain  
 on the complex manifold $X(\Bbb{C})$ with boundary on the divisor $D(\Bbb{C})$. Periods are the ring (over $\Bbb{Q}$) generated by numbers of the form
$\int_{\gamma}\omega$.
\end{definition}

We could have replaced $\Bbb{Q}$ by $\overline{\Bbb{Q}}$ above, and obtained the same ring (as Kontsevich and Zagier remark). This is easy because a variety
defined over $\overline{\Bbb{Q}}$ can be viewed as defined over $\Bbb{Q}$,
but we get several copies over the algebraic closure. But there is one more
modification that one can make which is a bit more subtle. This is to allow
for absolutely convergent integrals. Most examples (eg. multiple zeta values)
are not directly periods in the above sense, the integrals defining them can
have singularities on the boundary. To take care of this we note the following theorem which will be proved in Section 2.
\begin{theorem}\label{thm:periods3}
  Let $X$ be a smooth  
  $n$-dimensional algebraic variety defined over a 
  field $k\subset\Ralg$.  Let $F$ be a reduced effective divisor and let 
  $\omega\in\Omega^n(X-F)$ be an $n$-form.  Let $C\subset X(\RR)$ be a 
  pre-oriented
  semi-arithmetic set with non-empty interior $C^{o}$. 
  Then the integral $\int_C\omega\in\periods$ 
  provided that it is absolutely convergent.
\end{theorem}
\begin{remark} Already known to Kontsevich and Zagier, as in page 31 of ~\cite{KontZagier}. We wanted to elaborate on their comment that this follows from resolution of singularities in characteristic $0$.
\end{remark}

We now turn to the theorem on Igusa Zeta functions. 

\begin{theorem}
 \label{thm:Per} Let $X$ be a smooth variety defined over $k\subset\Ralg$ 
 and let $f\in\coords(X)$ be a function.  Let $C$
 be a compact pre-oriented semi-arithmetic subset of $X_{f \geq
 0}(\RR)$ defined over $k$.  Then, if $\omega\in\Omega^n(X)$ is a
 differential form, the function
\begin{equation}
  I(s)=\int_{C} f^s\omega 
\end{equation}
extends meromorphically to all of $\CC$ with poles occurring only at negative integers.
Moreover, for any $s_0\in\ZZ$, the coefficients $a_i$ in the Laurent expansion
\begin{equation}
  I(s)=\sum_{i\geq N} a_i (s-s_0)^i
\end{equation}
are periods.
\end{theorem}

Our first step is to prove the theorem for $s_0>0$.  In this case,
Atiyah's theorem shows that the integral for $I(s)$ converges and is
analytic in a neighborhood of $s_0$.  Thus, assuming $f\neq 0$, we can
differentiate under the integral sign to obtain
\begin{equation}
  \label{eq:OneC}
  I^{(l)}(s_0) = \int_{C} f^{s_0} \log^l (f)\omega.
\end{equation}

Now 
\begin{equation}
  \label{eq:OneD}
  \log f(x) = \int_0^1 \frac{f(x)-1}{(f(x)-1)t + 1}\, dt.
\end{equation}  
Thus we can write the $\log$ factors in (\ref{eq:OneC}) as period integrals.

To do this explicitly,  set $Y=X\times\Aff^l$,  $D=C\times [0,1]^l$ and 
$$
\eta=\omega\wedge\frac{f(x)-1}{(f(x)-1)t_1) + 1}\, dt_1\wedge\cdots
           \wedge\frac{f(x)-1}{(f(x)-1)t_l) + 1}\, dt_l.
$$
We then have 
\begin{equation}
  \label{eq:OneE}
\int_{D} f^{s_0}\eta =\int_{C} f^{s_0} \log^l (f)\omega.
\end{equation}

The left hand side is absolutely convergent (in fact bounded on the domain of integration\footnote{ Consider the integral $\int_{1\geq x\geq y\geq 0}\frac{y}{x}dxdy$, the integrand is bounded, yet the differential form $\frac{y}{x}dxdy$
has a pole at $x=0$. A blow up at $(0,0)$ resolves this problem and converts this to $\int_{0\leq x\leq 1, 0\leq u \leq 1} ux dx du$.}). Thus, $I^{(l)} (s_0)$ is a period for all $l$ as long as $s_0>0$ (theorem ~\ref{thm:periods3}), and the 
theorem is verified for $s_0>0$. 

To verify the theorem for $s_0\leq 0$, we use an auxilliary function and 
the Picard-Fuchs equation.  Set 
\begin{equation}
  \label{eq:Jdef}
J(t) = \int_C \frac{\omega}{1-tf}
\end{equation}
viewing the integrand as an $n$-form on $X\times\Aff^1$.   
Then $J(t)=\sum_{l\geq 0} I(l) t^l$ for all $t$ such that the sum converges.
Since $C$ is compact, $f$ is bounded on $C$ by some constant $R$.  Thus, 
for $t<1/R$, $C$ does not intersect the divisor $Z=V(1-tf)$ where the 
integrand may have a pole, and the integral (\ref{eq:Jdef}) converges.

Using the triangulation theorem for semi-algebraic sets
(\cite{bochnak} Theorem 9.2.1), we can assume that $C$ is homeomorphic
analytically to an $n$-simplicial complex with one $n$-cell and that
$\partial C$ is contained in a divisor $D\subset X$ (defined over $k$)
Let $\sigma\in
H_n(X(\RR)-Z(\RR),D(\RR) -Z(\RR);\ZZ)$ 
be the class represented by integration over the
points of $C$ that are smooth in $X$.  Then, for each $t$ with
$|t|<1/R$, 
\begin{equation}
  \label{eq:Jdefinition}
J(t)=\int_{\sigma} \frac{\omega}{1-tf}.
\end{equation}

There is an algebraic vector bundle $\mathcal{V}=H_{DR}^{m}(X-Z,D-Z)$ over $\Aff^1- S$ where $S$ is   a finite subset of $\Aff^1$ which can include 0 (but defined over $k$). The stalks of $\mathcal{V}$
are equal to the de Rham cohomology of $H^n(X_y-Z_y,Z_y-D_y)$ over the field
$k(y)$ for all $y\in \Aff^1-S$. The integrand $s=\frac{\omega}{1-tf}$ can be thought of as a global section of $\mathcal{V}$ (because it is an algebraic differential form of the top degree it is closed and vanishes when restricted to $D-Z$)

This bundle $\mathcal{V}$ carries an algebraic connection $\nabla$, an isomorphism
over $\Aff^1-S$ (of analytic vector bundles)
$$\mathcal{V}_{\Aff^1_{\Bbb{C}}}\to \mathcal{L}\bigotimes_{\Bbb{Z}} O_{Y_{\Bbb{C}}}$$ where
$\mathcal{L}$ is the local system whose fiber at $y\in  Y_{\Bbb{C}}$ is the singular cohomology of the pair $(X_y-Z_y,D_y-Z_y)$. The connection is integrable, has regular singular points and the sheaf of flat sections is the sheaf $\mathcal{L}$.

 If $\sigma$ is a flat section of the dual local system $\mathcal{L^*}$ (which is the local system of the homology of pairs $H^n(X_y-Z_y,D_y-Z_y)$) over an open set $U\subset \Aff^1-S$ in the analytic topology, then we can form a function on $U$: $g(y)=
\int_{\sigma} s_y$. If $T$ is a tangent vector field on $U$, we have the formula
$$T(g)=\int_{\sigma} \nabla_{T}(s)_y.$$

Now, $\mathcal{V}$ is a vector bundle of finite rank so given any section $s$ over $\Aff^1-S$, there is a relation of the form

$$\sum_{i=0}^r q_i(t)  \nabla^i_{T}(s)_y.$$

where the $q_i$ are rational functions in $t$ with coefficients in $k$. We can assume that they are polynomials by multiplying the equation by a polynomial $\in k[t]$.

Integrating this against the $\sigma$ obtained from $C$ and $s=\frac{\omega}{1-tf}$  we obtain a 
nontrivial linear relation of the form
\begin{equation}
  \label{eq:P-F}
  \sum_{i=0}^r q_i(t) J^{(i)} (t)=0
\end{equation}
where the $q_i(t)\in k(t)$.  

For a complete
reference to the Picard-Fuchs theory see~\cite{deligne-diffeq}.

Clearing denominators in (\ref{eq:P-F}), we can assume that 
the $q_i(t)\in k[t]$.   Expanding out 
$q_i(t)=\sum_{j=0}^{d_i} a_{i,j} t^j$  (for some $a_{i,j}\in k$)
and 
$\displaystyle J^{(i)}(t) = \sum_{j\geq 0} 
\frac{j!}{(j-i)!} t^{j-i} I(j)$ 
and equating  
terms with the same power of $t$, we obtain a
set of relations between the $I(j)'s$.  
Explicitly, we obtain the relation
\begin{equation}
  \label{eq:BigR}
  \sum_{s\geq 0}\sum_{i=1}^r\sum_{j=0}^{d}
  a_{i,j}\frac{(s+i-j)!}{(s-j)!}  I(s+i-j)\, t^s = 0
\end{equation}
where $d=\max{d_i}$.

Noting that, for each pair $(i,j)$, the coefficient 
$\displaystyle a_{i,j}\frac{(s+i-j)!}{(s-j)!} $
is a polynomial
of degree $i$ in $s$, we see that we have a relation of the form
\begin{equation}
  \label{eq:rel2}
  \sum_{i=0}^{d+r} c_i(s) I(s+i) =0
\end{equation}
with the $c_i$ polynomials in $k[s]$.  Note that, as long as $f$ and
$\omega$ are nonzero the relation (\ref{eq:rel2}) is nontrivial.

We wish to show tht (\ref{eq:rel2}) holds for all complex values of $s$.
By the uniquness of analytic continuation, it is enough to show that this is
so for $\Re(s)>d+r$.  We then use the following corollary of a 
result from~\cite{carleson} (p. 953).
\begin{theorem}[Carleson]
  Let $h(z)$ be holomorphic for $\Re(z)>0$ and assume $h(n)=0$ for 
  $n\in\NN$.  Then $h(z)=0$ if $h(z)\leq Ke^{m\Re(z)}$ for a constants $m$ and
  $K$.
\end{theorem}

To use Carleson's theorem, let $Q(z)$ be the left hand side of
(\ref{eq:rel2}) viewed as a function of a complex variable $z=s-d-r$.
Then $Q(z)$ is holomorphic for $\Re(z)>0$.  Moreover, since $f$ is
bounded on the semi-algebraic set $C$ by a number $R$, $|I(s)|$ is
bounded by $AR^{\Re (s)}$ for some constant $A$.  Thus $Q(z)$ is
bounded by $Ke^{m\Re(s)}$ for some constants $m$ and $K$.  It follows
from Carleson's theorem that $Q(z)=0$ for $\Re(s)>0$.  Thus, by
uniqueness of analytic continuation, it follows that $Q(z)=0$ for all
$z$.

Without loss of generality, we can assume that $c_0(s)$ in (\ref{eq:rel2})
is nonzero.  Then we have a relation
\begin{equation}
  \label{eq:rel4}
  I(s)=\sum_{i=1}^{d+r} l_i(s) I(s+i).
\end{equation}
where $\displaystyle l_i=\frac{-c_i(s)}{c_0(s)}$.
Using (\ref{eq:rel4}), we can complete the proof of Theorem~\ref{thm:Per} 
by descending induction on $s_0$.   For $s_0>0$, the theorem is established.
Suppose then that the theorem is established for $s_0>M$.
We can then use the Laurent expansions for the terms on right hand side of 
(\ref{eq:rel4}) to write out the Laurent expansion for the left hand side.
Using the fact that the $l_i$ are rational function in $k(t)$ and using 
the Laurent expansions of $I(s)$ at $s_0>M$, it is easy to see that 
the theorem holds for $s=s_0$.
\begin{remark}Historical precedents to the `functional equation' in theorem 1.8 should be noted. They first appear in Bernstein's paper ~\cite{bernstein}. Using the theory of $\mathcal{D}$-modules, he shows that if the domain of of integration was all of $\Bbb{R}^n$ and the polynomial function $f$, satisfied a growth rate of the form
$$|f(X)|\geq C ||X||^A$$

for $A>0$ and $||(x_1,\dots,x_n)||=\sum x_i^2$, then functions of the type

$$H(s)=\int_{\Bbb{R}^n}f^{-s}dx_1,\dots,dx_n$$ satisfied functional equations.
This was achieved beautifully using the theory of $\mathcal{D}-$modules. But this approach fails (or atleast we could not make it work) when the domain of integration is an arbitrary semi-algebraic set. 
\end{remark}

\section{Periods and semi-arithmetic sets}

In this section we prove a theorem tacitly used in~\cite{KontZagier}
relating integrals over semi-arithmetic sets to the integrals over
cohomology classes which are more widely thought of as period
integrals.  The main tool is the same corollary of resolution of
singularities used by Atiyah to prove theorem~\ref{thm:Atiyah}.  We
state it here in the form that we will use.
\begin{theorem}[Resolution Theorem]
  \label{thm:res}
Let $F\in\coords(X)$ be a nonzero function
on a smooth, complex $n$-dimensional algebraic variety.
Let $\omega\in\Omega^n(X-E)$ be a differential $n$-form where 
$E$ is a divisor.  Let $Z(\omega)$ denote the zero set of $\omega$.
Then there is a proper morphism $\varphi:\tilde{X}\to X$ 
from a smooth variety $\tilde{X}$ such that
\begin{enumerate} 
\item $\varphi: \tilde{X}-\tilde{A}\to X-A$ is an isomorphism,
where $A=F^{-1}(0)\cup E\cup Z(\omega)$ and $\tilde{A}=\varphi^{-1}(A)$.
\item for each $P\in\tilde{X}$ there are local coordinates 
$(y_1,\ldots, y_n)$
centered at $P$ so that, locally near $P$, 
\begin{eqnarray*}
F\circ\varphi &=& \epsilon\prod_{j=1}^n y_j^{k_j}\\
\omega        &=& \delta\prod_{j=1}^n y_j^{l_j}\, 
                   dy_1\wedge\cdots\wedge dy_n 
\end{eqnarray*}
\end{enumerate}
where $\epsilon,\delta$ are units in $\coords_{X,P}$, the $k_j$ are 
non-negative integers and the $l_j$ are arbitrary integers.
\end{theorem}

The theorem, the statement of which is very close to the statement of
Atiyah's resolution theorem on p. 147 of~\cite{Atiyah70}, is proved by
applying Main Theorem II in ~\cite{hironaka} to the ideals $F\coords_X$,
$E$ and $Z(\omega)$.

\begin{proposition} Let $X$ be a smooth $n$-dimensional algebraic
variety defined over $\Ralg$.  Let $F$ be a reduced effective divisor
and let $\omega\in\Omega^n(X-E)$ be an $n$-from.  Let 
$$
G=\{x\in X(\RR) | g_i(x)\geq 0\}
$$
for some set $\{g_i\}_{i=1}^m$ of functions in $\coords(X)$.
be a compact, pre-oriented semi-algebraic set with non-empty interior
$G^{0}$.  Then $\int_G\omega$ converges absolutely only if there is 
a smooth $n$-dimensional algebraic variety 
$\tilde{X}$ with proper, birational morphism $\varphi:\tilde{X}\to X$ and 
a compact semi-algebraic set $\tilde{G}$ such that 
\begin{enumerate}
\item $\int_{\tilde{G}} \varphi^*\omega = \int_G \omega$.
\item $\varphi^*\omega$ is holomorphic on $\tilde{G}$.
\end{enumerate}
\end{proposition}
\begin{proof}  
Using the resolution theorem with $F=\prod_{i=1}^m g_i$, 
we can find a smooth variety
$\tilde{X}$ with a proper, birational morphism to $X$ such that 
for every point $P\in\tilde{X}$ we have local parameters 
$(y_1,\cdots, y_n)$ defined in a neighborhood of $P$
with
\begin{eqnarray*}
g_i\circ\varphi &=& \epsilon_i \prod_{j=1}^n y_j^{k_{ij}}\\
\varphi^*\omega &=& \delta     \prod_{j=1}^n y_j^{l_j}\, 
                                          dy_1\wedge\cdots\wedge dy_n.
\end{eqnarray*}
Here the $\epsilon_i$ and $\delta$ are invertible near $P$.
Set $\tilde{G}$ equal to the analytic closure of $\varphi^{-1}(G-A)$
with $A$ as in the resolution theorem.  Then 
$\int_{\tilde{G}}\varphi^*\omega=\int_G\omega$ because $\tilde{G}$ and 
$G$ differ only by measure $0$ sets.  Moreover, since $\varphi$ is 
proper and $\tilde{G}$ is a closed subset of  $\varphi^{-1} G$,
$\tilde{G}$ is compact.

To see that $\varphi^*\omega$ is holomorphic on $\tilde{G}$, let 
$P\in\tilde{G}$ be a point and let $\tilde{U}$ be a neighborhood
of $P$ with a local coordinate system $(y_1,\cdots, y_n)$ 
as in the resoltuion theorem.  Since $P$ is in the closure of 
$\varphi^{-1}(G-A)$, $g_i(P)\geq 0$ for all $i$.  Let 
$s_i$ be the sign ($\pm 1$) of $\epsilon_i(P)$.  Then, since
$\int_{G}\omega$ is absolutely convergent, it follows that 
\begin{equation}
\int_{0<s_iy_i(p)<r} \varphi^*\omega=
\int_{0<s_iy_i(p)<r} \prod_{j=1}^n y_j^{l_j}\, 
                                dy_1\wedge\cdots\wedge dy_n
\end{equation}
is absolutely convergent for a sufficiently small $r$.
It is easy to see that this is not possible unless $l_j\geq 0$ for
all $j$.  Thus $\varphi^*\omega$ is holomorphic at $P$.
\end{proof}

\begin{proposition} Let $X$ be a smooth algebraic variety over 
$\Ralg$ and let $G=\{x\in X(\RR) | g_i(x)\geq 0\}$ be a compact 
pre-oriented set.  
Let $\omega\in\coords_X(X)$ be a differential $n$-form.  
Then there is a divisor $D\subset X$ and a chain 
$\sigma\in H_n(X,D)$ such that $\int_G\omega=\int_{\sigma} \omega$.
\end{proposition}
\begin{proof}
The pre-orientation on $G$ gives us a dense, smooth, open semi-algebraic 
subset $U$ in $G$ with an orientation on $U$.  We, therefore, obtain
a chain $\sigma\in H_n(X,D)$ where $D$ is the set of zeroes of the 
functions $g_i$ defining $G$.  This $\sigma$ corresponds to the orientation on the 
open subset $U$ so we have $\int_{\sigma} \omega = \int_G\omega$.
\end{proof}
Therefore we conclude:
\begin{theorem}\label{thm:periods}
  Let $X$ be a smooth  
  $n$-dimensional algebraic variety defined over a 
  field $k\subset\Ralg$.  Let $F$ be a reduced effective divisor and let 
  $\omega\in\Omega^n(X-F)$ be an $n$-form.  Let $C\subset X(\RR)$ be a 
  pre-oriented
  semi-arithmetic set with non-empty interior $C^{o}$. 
  Then the integral $\int_C\omega\in\periods$ 
  provided that it is absolutely convergent.
\end{theorem}

\bibliographystyle{plain}

\begin{thebibliography}{10}

\bibitem{Atiyah70}
 M. F. Atiyah, M. F.
 \newblock Resolution of singularities and division of distributions.
 \newblock Comm. Pure Appl. Math,23,1970, 145--150.

\bibitem{BernsteinGelfand69}
 Bernstein, I. N. and Gelfand, S. I.
 \newblock Meromorphy of the function ${P}\sp{\lambda }$.
 \newblock Funkcional. Anal. i Prilov zen.,3,1969,1,84--85.

 \bibitem{bernstein}
 Bern{\v{s}}te{\u\i}n, I. N.
 \newblock Analytic continuation of Generalised functions.
 \newblock Functional Anal. Appl,6, 1972, 273--285(1973).

\bibitem{bochnak}
 J. Bochnak and M. Coste and Marie-Fran{\c{c}}oise Roy.
 \newblock Real algebraic geometry.
 \newblock Ergebnisse der Mathematik und ihrer Grenzgebiete (3) [Results
              in Mathematics and Related Areas (3)].
 \newblock volume  {36}, Springer-Verlag, Berlin 1998.

\bibitem{carleson}
 L. Carleson.
 \newblock On Bernstein's approximation problem.
 \newblock Proc. Amer. Math. Soc.,2, 1951, 953--961.

\bibitem{deligne-diffeq}
 Pierre Deligne.
 \newblock \'{E}quations diff\'erentielles \`a points singuliers
              r\'eguliers.
 \newblock Lecture Notes in Mathematics, Vol. 163, Springer-Verlag, Berlin 1970.

\bibitem{hironaka}
H. Hironaka.
\newblock Resolution of Singularities of an algebraic variety over a field of characteristic 0.
\newblock Annals of Math, 79, 1964, 109--326.

\bibitem{KontZagier}
M. Kontsevich and D. Zagier.
\newblock Periods.
\newblock Mathematics Unlimited --- 2001 and Beyond 771 -- 809, Springer-Verlag, Berlin, 2001.





\end{thebibliography}

\def\noopsort#1{}

\end{document}